\RequirePackage{amsmath}
\documentclass[dvipsnames,runningheads]{llncs}
\usepackage{bm}
\usepackage[T1]{fontenc}
\usepackage{lmodern}
\usepackage{graphicx}
\usepackage{commands}
\usepackage[hidelinks]{hyperref}
\usepackage{color}
\usepackage[export]{adjustbox}

\begin{document}
\title{Diffusion-Shock Filtering on the Space of Positions and Orientations}
\author{
Finn M. Sherry\inst{1}
\and
Kristina Schaefer\inst{2}
\and 
Remco Duits\inst{1}
}
\authorrunning{F.M. Sherry et al.}
\institute{CASA \& EAISI, Dept. of Mathematics \& Computer Science, Eindhoven University of Technology, the Netherlands\\
\email{\{f.m.sherry,r.duits\}@tue.nl}
\and
Mathematical Image Analysis Group, Dept. of Mathematics \& Computer Science,
E1.7, Saarland University, 66123 Saarbrücken, Germany\\
\email{schaefer@mia.uni-saarland.de}
}
\maketitle              

\begin{abstract}
We extend Regularised Diffusion-Shock (RDS) filtering from Euclidean space $\Rtwo$ \cite{schaefer2024regularisedds} to the space of positions and orientations $\M \coloneqq \Rtwo \times S^1$.
This has numerous advantages, e.g. making it possible to enhance and inpaint crossing structures, since they become disentangled when lifted to $\M$. 
We create a version of the algorithm using gauge frames to mitigate issues caused by lifting to a finite number of orientations. This leads us to study generalisations of diffusion, since the gauge frame diffusion is not generated by the Laplace-Beltrami operator. 
RDS filtering compares favourably to existing techniques such as Total Roto-Translational Variation (TR-TV) flow \cite{smets2021tvflow,chambolle2019tv}, NLM \cite{nlm}, and BM3D \cite{bm3dold} when denoising images with crossing structures, particularly if they are segmented. Additionally, we see that $\M$ RDS inpainting is indeed able to restore crossing structures, unlike $\Rtwo$ RDS inpainting.
\end{abstract}

\section{Introduction}
With the intention of creating a transparent inpainting model with stability guarantees, Schaefer and Weickert introduced Diffusion-Shock (DS) filtering in~\cite{SW23} and refined it to Regularised Diffusion-Shock (RDS) filtering in~\cite{schaefer2024regularisedds}. It combines homogeneous diffusion~\cite{Ii62} with coherence enhancing shock filtering~\cite{We03}. While RDS filtering was introduced primarily for inpainting, it can also be used for denoising; the presence of an explicit shock term makes it capable of edge-preserving denoising. In this article, we extend RDS filtering from $\Rtwo$ to the space of positions and orientations $\M \coloneqq \Rtwo \times S^1$ to combine the edge-preserving denoising capabilities of RDS filtering with the benefits of processing with partial differential equations (PDEs) on $\M$, such as preservation of crossings \cite{franken2009cedos,smets2021tvflow}.

Many PDE based image processing techniques have been extended to spaces of positions and orientations $\mathbb{M}_d \coloneqq \mathbb{R}^d \times S^{d - 1}$.
Total Roto-Translational Variation (TR-TV) on $\M$ has been studied by Chambolle \& Pock \cite{chambolle2019tv} and Smets et al.~\cite{smets2021tvflow}; see Pragliola et al.~\cite{CalatroniSIAMreview2023} for a review on the topic. There are interesting links between TV flow and elastica which have inspired further regularisation methods \cite{Liu2023ElasticaRegularization,chambolle2019tv}.
Mean Curvature (MC) flows on $\mathbb{M}_d$ have been employed for various purposes including denoising and inpainting by Citti et al. \cite{Citti2016mcf} ($d = 2$) and St.~Onge et al. \cite{stonge2019mcf} ($d = 3$). 
Diffusion PDEs on $\mathbb{M}_d$ have been employed in denoising, inpainting, and neurogeometry \cite{cittisartiJMIV,Prandi,franken2009cedos,petitotbook}.
Morphology HJB-PDEs on $\mathbb{M}_d$ have been effective in geometric deep learning \cite{smets2024thesis,paperGijs} and fiber enhancement \cite{DuitsJMIV2013}.
RDS filtering combines diffusion for denoising with morphological shock filtering to preserve edges like TR-TV \cite{chambolle2019tv} and MC flows \cite{smets2021tvflow}; by extending the RDS processing to $\M$ we include preservation of crossing/overlapping structures.

\subsubsection{Our Contributions.}
Two novel methods based on $\Rtwo$ RDS filtering are introduced.
We create crossing aware RDS filtering by applying it in the space of positions and orientations $\M$.
This gives rise to our new evolution operators:
\mbox{}~\hspace{0.5 cm} (1) left-invariant RDS filtering, and \\
\mbox{}~\hspace{0.5 cm} (2) gauge frame RDS filtering. \\ 
We analyse our evolutions theoretically. In the left-invariant setting, the Laplacian used for the diffusion part of the RDS filtering is both the Lie-Cartan Laplacian and the Laplace-Beltrami operator, cf.~Theorem~\ref{thm:left_invariant_lie_cartan_laplacian}.
However, for the gauge frame RDS filtering, a discrepancy occurs, which has been overlooked in \cite{franken2009cedos,smets2021tvflow}. We derive the difference between these Laplacians in Theorem~\ref{thm:data_driven_lie_cartan_laplacian}.

We verify the behavior of our novel methods experimentally. Our denoising experiments show that our novel methods outperform other PDE-based methods, like $\Rtwo$ RDS filtering, TR-TV flows on $\M$~\cite{smets2021tvflow} as well as non-PDE-based denoising methods \cite{bm3dold,nlm} in terms of Peak Signal-to-Noise Ratio (PSNR) in denoising tasks. Moreover we confirm that  $\M$ RDS filtering allows for inpainting crossing structures unlike its $\Rtwo$ version.
The implementations of both methods and experiments are available at \url{https://github.com/finnsherry/M2RDSFiltering}.
\section{Regularised Diffusion-Shock Filtering on \texorpdfstring{$\Rtwo$}{Euclidean Space}}\label{sec:r2_ds}
RDS filtering combines homogeneous diffusion with coherence-enhancing shock filtering. 
The coherence-enhancing shock filter sharpens and elongates edge-like structures by adaptively applying the morphological operations \emph{dilation} and \emph{erosion}~\cite{So04}. Dilation of a greyscale image $f: \Omega \subset \Rtwo \to \R$ replaces the grey value in location $\vec{x}$ by the supremum of $f$ within a specified neighbourhood around $\vec{x}$. Erosion uses the infimum instead. The PDE-based formulation~\cite{BM92,AGLM93,AVK93} of dilation $(+)$ / erosion ($-$) with a disk-shaped neighbourhood is given by
\begin{equation}
\partial_t u = \pm \abs{\gradient u},
\end{equation}
with the initial image $u(\vec{x}, 0) = f(\vec{x})$ and reflecting boundaries, where $\abs{\,\cdot\,}$ and $\gradient$ are the Euclidean norm and gradient, respectively. 
The coherence-enhancing shock filter applies dilation when the data is concave in the direction perpendicular to the local orientation, and erosion when it is convex. 
This direction is determined by the dominant eigenvector $\vec{w}$ (i.e. the eigenvector with the largest eigenvalue) of a structure tensor~\cite{FG87} $\mat{J}_\rho(\gradient u_\sigma) = K_\rho * (\gradient u_\sigma \gradient u_\sigma^\top)$ where $u_\sigma = K_\sigma *u$ with the Gaussian convolution kernels $K_\rho, K_\sigma$. With that, the coherence-enhancing shock filter evolves an initial greyscale image $f$ by
\begin{equation}\label{eq:shock_R2}
\partial_t u = -S(\partial_{\vec{w}\vec{w}} u_\sigma) \abs{\gradient u},
\end{equation}
with initial condition $u(\vec{x}, 0) = f(\vec{x})$ and reflecting boundaries.
The sigmoidal function $S: \R \to [-1, 1]$ implements the behaviour of a (soft) sign function. 

RDS filtering aims at applying this shock filter near edges, while the diffusion smooths flat areas. This adaptive behaviour is produced using a Charbonnier weight~\cite{charbonnier1997switch} $g: \R_{\geq 0} \to \R: x \mapsto \sqrt{1 + x/\lambda^2}^{-1}$ with a Gaussian-smoothed gradient magnitude $\gradient u_\nu$ as input. 
In summary, the RDS filtering PDE is given by
\begin{equation}\label{eq:ds_R2}
\partial_t u = g\left(\abs{\gradient u_\nu}^2\right) \laplace u
  - \Big(1 - g\left(\abs{\gradient u_\nu}^2\right)\Big)
  S\left(\partial_{\vec{w}\vec{w}} (u_\sigma) \right)
  \abs{\gradient u},
\end{equation}
with initial condition $u(\vec{x}, 0) = f(\vec{x})$ and reflecting boundaries. 

For the application to digital images, the PDE can be discretised with an explicit scheme, that inherits a maximum-minimum principle which excludes under/over shoots. Diffusion is discretised with central differences, the morphological terms requires an upwind scheme \cite{rouy1992viscosity}. For details see~\cite{SW23,schaefer2024regularisedds}. 
\section{Space of Positions \& Orientations\texorpdfstring{ $\M$}{}}\label{sec:m2}
\begin{figure}
\centering
\begin{subfigure}[t]{0.2\textwidth}
\begin{tikzpicture}
\draw(0, 0)node[inner sep=0]{\includegraphics[width=\textwidth]{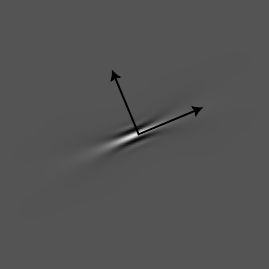}};
\draw(0.75, 0.23)node{$\A_1$};
\draw(-0.18, 0.7)node{$\A_2$};
\end{tikzpicture}
\caption{Cake wavelet.}\label{fig:cakewavelet}
\end{subfigure}
~ 
\begin{subfigure}[t]{0.75\textwidth}
\begin{tikzpicture}
\draw(0, 0)node[inner sep=0]{\includegraphics[width=\textwidth]{schematic_road_block_combined_horizontal.pdf}};
\draw(-3.81, -0.1)node{$\mathcal{W}_\psi$};
\draw(-4.05, 0.77)node{$\theta$};
\draw(0.0, 0.6)node{Processing};
\draw(0.67, -0.1)node{$\mathcal{W}_\psi^{-1}$};
\draw(4.23, 0.77)node{$\theta$};
\end{tikzpicture}
\caption{Performing multi-orientation processing.}\label{fig:disentanglement}
\end{subfigure}
\caption{(a) Cake wavelet for orientation $\theta = \frac{\pi}{8}$. (b) Lifting disentangles crossing and overlapping structures.}
\label{fig:summary_M2}
\end{figure}

Multi-orientation image processing has been beneficial in denoising \cite{chambolle2019tv,stonge2019mcf} and line enhancement \cite{Citti2016mcf,franken2009cedos}. We lift the image data to the \emph{space of positions and orientations}, $\M$, with an \emph{orientation score transform}, $\mathcal{W}_\psi$, which disentangles crossing and overlapping structures. This allows for orientation-aware image processing in $\M$. After processing, the result is mapped back to $\Rtwo$ with an inverse transform $\mathcal{W}_\psi^{-1}$. Figure~\ref{fig:disentanglement} illustrates this process. 
\begin{definition}[Space of positions and orientations\texorpdfstring{ $\M$}{}]\label{def:m2}
The \emph{space of positions and orientations} is defined as the smooth manifold $\M \coloneqq \Rtwo \times S^1$. 
Elements of $\M$ are denoted by $(\vec{x}, \theta) \in \Rtwo \times \R / 2\pi \Z$.
\end{definition}
The Lie group of roto-translations $\SE(2) \coloneqq \Rtwo \rtimes \SO(2)$ acts on $\M$ as follows: for $g = (\vec{x}_g, R_{\theta_g}) \in \SE(2)$ and $p = (\vec{x}, \theta) \in \M$ we have action $L_g p \coloneqq (R_{\theta_g}\vec{x} + \vec{x}_g, \theta + \theta_g)$, where $R_\theta$ is a counter-clockwise rotation by $\theta$. By choosing reference element $p_0 \coloneqq (\vec{0}, 0) \in \M$, we see that \cite{smets2021tvflow}
\begin{equation}\label{eq:m2_principal_homogeneous_space}
\M \cong \SE(2) / \Stab_{\SE(2)}(p_0) \cong \SE(2) / \{e\} \cong \SE(2),
\end{equation}
with unit element $e \coloneqq (\vec{0}, I)$, so $\M$ is the \emph{principal homogeneous space} of $\SE(2)$. 
On $\M$, we perform $\SE(2)$-equivariant processing via \emph{left-invariant vector fields},
\begin{equation}\label{eq:left_invariant_vector_fields}
\mathfrak{X}(\M) \coloneqq \{\A \in \sections(T\M) \mid (L_g)_* \A_p = \A_{L_g p}, \forall g \in \SE(2), p \in \M\},
\end{equation}
where $\sections(T\M)$ are the smooth sections of the tangent bundle $T\M$, i.e. the set of smooth vector fields on $\M$.
\begin{definition}[Left-Invariant Frame]\label{def:left_invariant_frame}
We define the left-invariant vector fields $\A_1, \A_2, \A_3 \in \mathfrak{X}(\M)$ as 
\begin{equation}\label{eq:left_invariant_vector_frame}
\A_1|_\mathbf{p} \coloneqq (L_{g_\mathbf{p}})_* \partial_x|_\mathbf{p_0}, \A_2|_\mathbf{p} \coloneqq (L_{g_\mathbf{p}})_* \partial_y|_\mathbf{p_0}, \textrm{ and } \A_3|_\mathbf{p} \coloneqq (L_{g_\mathbf{p}})_* \partial_\theta|_\mathbf{p_0},
\end{equation}
where $g_\mathbf{p} = (\mathbf{x}, R_\theta)$ for $\mathbf{p} = (\mathbf{x}, \theta)$, and $\cdot_*$ denotes the pushforward.
Together, they form a basis for $\mathfrak{X}(\M)$, which we call the \emph{left-invariant frame}.
\end{definition}
We can gain the benefits of $\M$ processing for data on $\Rtwo$ by lifting the data to $\M$, using the \emph{orientation score transform}.
\begin{definition}[Orientation Score]\label{def:orientation_score}
The \emph{orientation score transform} $\mathcal{W}_\psi : \Ltwo(\Rtwo) \to \Ltwo(\M)$, where $\psi$ is a \emph{cake wavelet} (see Figure~\ref{fig:cakewavelet}), is defined by
\begin{equation}\label{eq:ost}
\mathcal{W}_\psi f (\vec{x}, \theta) \coloneqq \int_{\Rtwo} \overline{\psi(R_{\theta}^{-1} (\vec{y} - \vec{x}))} f(\vec{y}) \diff \vec{y}
\end{equation}
for $f \in \Ltwo(\Rtwo)$ and $(\vec{x}, \theta) \in \M$. We then call $\mathcal{W}_\psi f$ the \emph{orientation score} of $f$.
\end{definition}
By design \cite{franken2009cedos}, cake wavelets allow for fast approximate reconstruction using
\begin{equation}\label{eq:fast_reconstruction}
f(\vec{x}) = \mathcal{W}_\psi^{-1}(\mathcal{W}_\psi f) (\vec{x})  \approx \int_{S^1} \mathcal{W}_\psi f (\vec{x}, \theta) \diff \theta.
\end{equation}
Working in the left-invariant frame has other upsides in addition to equivariance. By construction, $\A_1$ points spatially along the local orientation, while $\A_2$ points laterally in an orientation score made with cake wavelets \cite{franken2009cedos}. 
Additionally, by lifting to $\M$ crossings are disentangled, opening the door for e.g. inpainting crossing structures, which is impossible in $\Rtwo$, see Figure~\ref{fig:disentanglement}.
\section{Regularised Diffusion-Shock Filtering on \texorpdfstring{$\M$}{Space of Positions and Orientations}}\label{sec:m2_ds}
We will now go over how to extend RDS filtering to $\M$, starting by discussing how to perform diffusion on $\M$. We use Laplacians of the form
\begin{equation}\label{eq:left_invariant_diffusion_scheme_exact}
\laplace_{\mathcal{G}} = g^{11} \A_1^2 + g^{22} \A_2^2 + g^{33} \A_3^2,
\end{equation}
where $g^{ij}$ is the $i,j$-th component of some left-invariant dual metric tensor field $\mathcal{G}^{-1}$ with respect to the left-invariant frame, i.e. $\mathcal{G}$ is given by $\mathcal{G}(\mathcal{A}_i, \mathcal{A}_{j})= g_{ii} \delta_{ij} = |g^{ii}|^{-1}\delta_{ij}$. We discretise \eqref{eq:left_invariant_diffusion_scheme_exact} using second order central differences and linear interpolation.
We discretise the time derivative with a forward difference with a time step size $\tau$. The  resulting explicit scheme satisfies a maximum-minimum principle if $\tau^{-1} \geq 4 \left(\frac{g^{11} + g^{22}}{\Delta_{xy}^2} + \frac{g^{33}}{\Delta_\theta^2}\right) \eqqcolon \tau_D^{-1}$, where $\Delta_{xy}$ and $\Delta_\theta$ are the spatial and orientational step sizes.

Let us now consider the shock filter. Like in $\Rtwo$, morphological dilations and erosions are generated by the norm of the gradient. Choosing a left-invariant metric tensor field $\mathcal{G}$ that is diagonal w.r.t. the left-invariant frame, we see that
\begin{equation}\label{eq:left_invariant_shock_scheme_exact}
\norm{\gradient_\mathcal{G} U}_\mathcal{G}^2 = g^{11} \abs{\A_1 U}^2 + g^{22} \abs{\A_2 U}^2 + g^{33} \abs{\A_3 U}^2.
\end{equation}
We discretise \eqref{eq:left_invariant_shock_scheme_exact} using a Rouy-Tourin-style upwind scheme \cite{rouy1992viscosity,bekkers2015subriemanniangeodesics}.
We discretise the time derivative with a forward difference with a time step size $\tau$. The  resulting explicit scheme satisfies a maximum-minimum principle if $\tau^{-1} \geq 2 \sqrt{\frac{g^{11} + g^{22}}{\Delta_{xy}^2} + \frac{g^{33}}{\Delta_\theta^2}} \eqqcolon \tau_S^{-1}$. 

The strength of RDS filtering in $\Rtwo$ comes from the guidance terms which allow the data to instruct whether to locally perform diffusion or shock, and dilation or erosion. As in the $\Rtwo$ case, we use the edge-preserving weight function by Charbonnier et al. \cite{charbonnier1997switch} to switch between diffusion and shock, 
\begin{equation}\label{eq:ds_switch}
g(\norm{\gradient_\mathcal{G} U}_\mathcal{G}^2) \coloneqq \sqrt{1 + \norm{\gradient_\mathcal{G} U}_\mathcal{G}^2 / \lambda^2}^{\, -1}.
\end{equation}
For the coherence-enhancing shock filter, we note that the local convexity can be determined by computing the Laplacian perpendicular to the local orientation: 
\begin{equation}\label{eq:morph_switch}
S(\laplace^\perp_{\mathcal{G}} U) \coloneqq S(g^{22} \A_2^2 U + g^{33} \A_3^2 U),
\end{equation}
with $S$ a sigmoidal function as in Equation~\eqref{eq:shock_R2}. Altogether, the evolution PDE is
\begin{equation}\label{eq:ds_SE22}
\partial_t U = g\left(\norm{{\bm \nabla_\mathcal{G}} U_\nu}_\mathcal{G}^2\right)\, \laplace_\mathcal{G} U
- \Big(1 - g \left(
   \norm{ {\bm \nabla_\mathcal{G}} U_\nu}_{\mathcal{G}}^2\right)\Big) 
 S_\rho \left(\laplace^\perp_\mathcal{G}U_\sigma \right)
  \norm{{\bm \nabla_\mathcal{G}} U }_{\mathcal{G}},
\end{equation}
where we have regularised the guidance terms, namely:
\begin{equation*}
U_\nu \coloneqq G_\nu *_{\SE(2)} U, \textrm{ and }
S_\rho(\laplace_\mathcal{G}^{\perp} U_\sigma) \coloneqq G_\rho *_{\SE(2)} S(\laplace_\mathcal{G}^{\perp} G_\sigma *_{\SE(2)} U),
\end{equation*}
with $G_\alpha$ a spatially isotropic Gaussian kernel with scale $\alpha > 0$ and $*_{\SE(2)}$ the group convolution on $\M$ (see \cite{franken2009cedos} for details on regularisation on $\SE(2) \cong \M$). The derivatives in the guidance terms are computed using central differences and linear interpolation. Finally, we apply reflecting spatial boundary conditions. The ranges of $g$ and $S$ are within $[-1, 1]$, so that the entire scheme is stable if the diffusion and the shock are stable individually; we hence choose $\tau \leq \min\{\tau_D, \tau_S\}$.

\subsubsection{Gauge Frames.}
In practice, we lift with a finite number of rotated cake wavelets, which means that the orientation scores and vector fields are discretised in the orientational direction. As such, data cannot always be lifted to exactly the correct orientation, and the vector field $\A_1$ only approximately points along the local orientation; the angle between the true spatial orientation and $\A_1$ is called \emph{deviation from horizontality} \cite{franken2009cedos}. 
Additionally, the direction of the lifted data will have an orientational component to account for the curvature in the input image, while $\A_1$ is purely spatial. 
Hence, it can be beneficial to use \emph{gauge frames}, which are adapted to the data. 
\begin{definition}[1st Gauge Vector]\label{def:first_gauge_vector}
Let $U \in C^2(\M)$, and let $\mathcal{G}$ be a metric tensor field on $\M$. The 1st gauge vector is given by (see \cite[Sec.~2.4]{smets2021tvflow} for details)
\begin{equation}\label{eq:first_gauge_vector}
\A_1^U|_p \coloneqq \underset{\substack{X_p \in T_p \M \\ \norm{X_p}_{\mathcal{G}} = 1}}{\argmin} \norm*{\nabla_{X_p}^{[0]} \gradient_{\mathcal{G}} U}_{\mathcal{G}}^2,
\end{equation}
where $\nabla_\cdot^{[0]}$ is the $0$ Lie-Cartan connection, which will be discussed in Sec.~\ref{sec:generalised_laplacians}.
\end{definition}
We choose $\mathcal{G} \coloneqq \mathcal{G}_\xi$, with $\mathcal{G}_\xi(\A_1, \A_1) = \xi^2 = \mathcal{G}_\xi(\A_2, \A_2)$ and $\mathcal{G}_\xi(\A_3, \A_3) = 1$, and $\mathcal{G}_\xi(\A_i, \A_j) = 0$ for $i \neq j$. $\xi$ has a large influence on the fitted gauge frame: we always choose $\xi = 0.1$. The other gauge vectors are then defined as follows: $\A_2^U$ is a purely spatial unit vector that is perpendicular to $\A_1^U$, and $\A_3^U$ is a unit vector perpendicular to both $\A_1^U$ and $\A_2^U$, s.t. $\{\A_i^U\}_i$ is a right-handed frame.

\begin{figure}
\centering
\begin{subfigure}[t]{0.45\textwidth}
\begin{tikzpicture}
\draw(0, 0)node[inner sep=0]{\includegraphics[width=\textwidth]{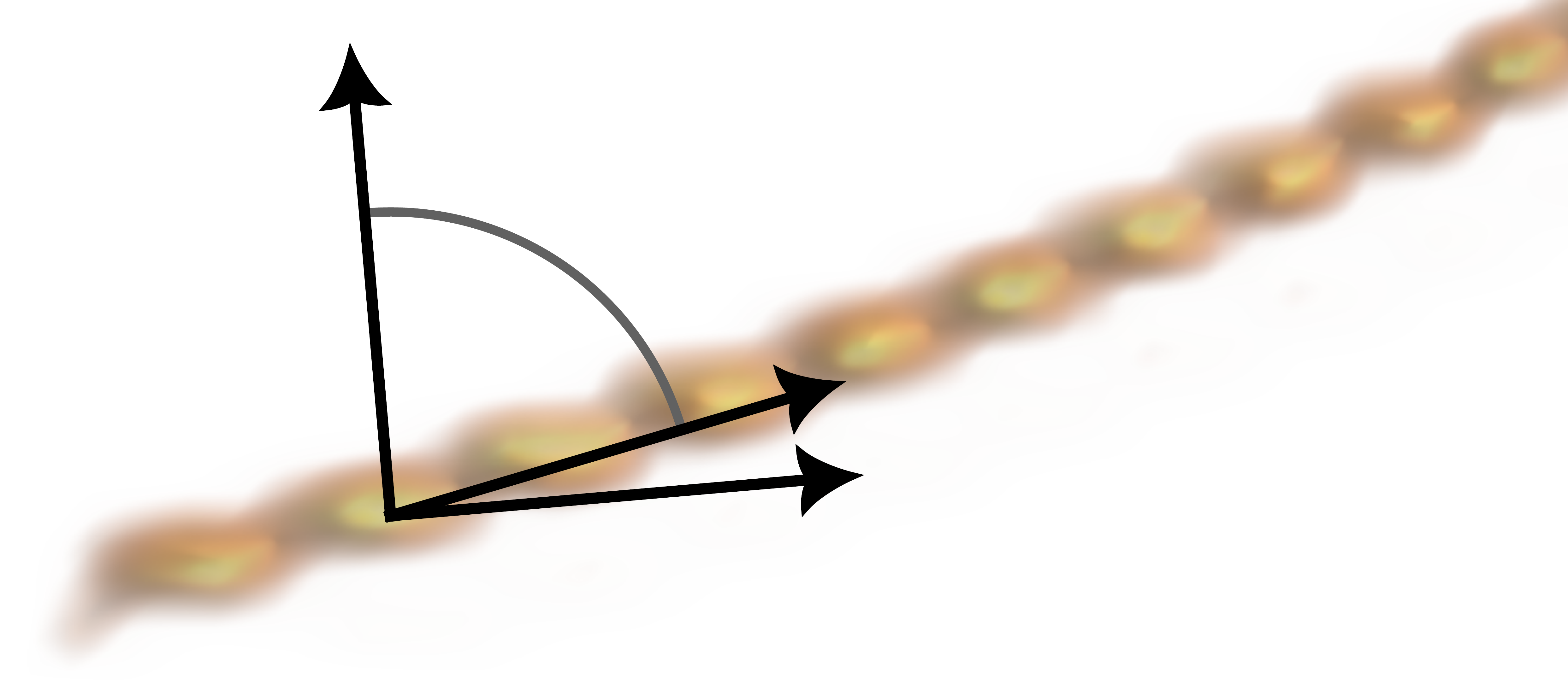}};
\draw(0.45, -0.5)node{$\A_1$};
\draw(0.42, -0.1)node{$\A_1^U$};
\draw(-1.45, 1.15)node{$\A_2$};
\draw(-0.15, 0.45)node{\textcolor{Maroon}{$\frac{\pi}{2} - d_H$}};
\end{tikzpicture}
\caption{Top View (along $\theta$-axis)}\label{fig:left_invariant_vs_gauge_top}
\end{subfigure}
~ 
\begin{subfigure}[t]{0.45\textwidth}
\begin{tikzpicture}
\draw(0, 0)node[inner sep=0]{\includegraphics[width=\textwidth]{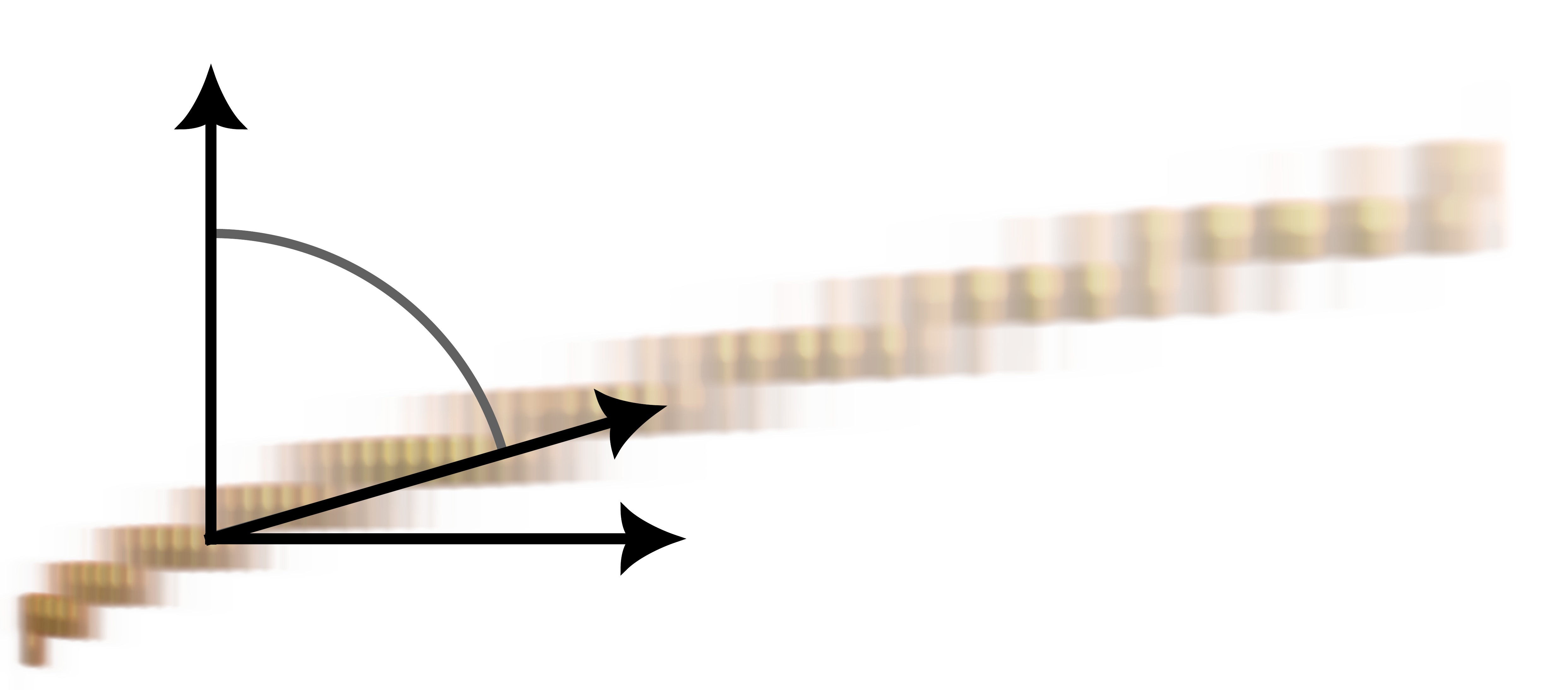}};
\draw(-0.14, -0.7)node{$\A_1$};
\draw(-0.17, -0.16)node{$\A_1^U$};
\draw(-1.95, 1.1)node{$\A_3$};
\draw(-0.3, 0.4)node{\textcolor{Maroon}{$\frac{\pi}{2} - \arctan(\kappa)$}};
\end{tikzpicture}
\caption{Side view (along $\mathcal{A}_{2}$-axis)}\label{fig:left_invariant_vs_gauge_side}
\end{subfigure}
\caption{Comparison of gauge frame and standard left-invariant frame. From the top view (a) we see that $\A_1^U$ has been rotated towards $\A_2$ to compensate for the deviation from horizontality $d_H$. From the side view (b) we see that $\A_1^U$ has been rotated towards $\A_3$; the rotation angle is related to the curvature $\kappa$.}
\label{fig:summary_gauge_frame}
\end{figure}
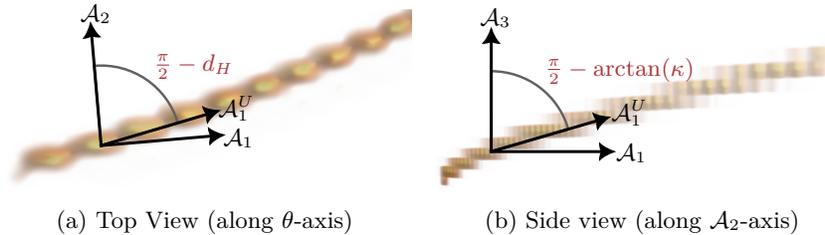

For the implementation of gauge frame RDS filtering, we can essentially copy the left-invariant PDE \eqref{eq:ds_SE22}, replacing the left-invariant frame vector fields $\A_i$ with the corresponding gauge frame vector fields $\A_i^U$. This scheme is again min-max stable if we use a timestep $\tau \leq \min\{\tau_D, \tau_S\}$.
\section{Generalised Laplacians}\label{sec:generalised_laplacians}
In Section~\ref{sec:m2_ds} we explained \emph{practically} how we generalise diffusion to $\M$; now we will investigate \emph{theoretically} how these diffusions compare to the typical notion of diffusion. 
On a Riemannian manifold $(M, \mathcal{G})$, one often defines diffusion as the evolution generated by the Laplace-Beltrami operator $\laplace_\mathcal{G} \coloneqq \divergence_\mathcal{G} \after \gradient_\mathcal{G}$. Motivated by the result that $\divergence_\mathcal{G}(X) = \trace(\nabla_{\cdot}^\mathrm{LC} X)$, where $\nabla^\mathrm{LC}$ is the Levi-Civita connection corresponding to $\mathcal{G}$ and $X \in \sections(T M)$ \cite[Ch.~5]{lee2018riemannian}, we define Laplace operators using general affine connections $\nabla$: 
\begin{definition}[Generalised Laplacian]\label{def:generalised_laplacian}
Let $(M, \mathcal{G})$ be a Riemannian manifold, and $\nabla$ an affine connection thereon. Then we define the corresponding \emph{generalised Laplacian} as
\begin{equation}\label{eq:general_connection_laplacian}
\laplace_{\mathcal{G}, \nabla} \coloneqq \divergence_\nabla \after \gradient_\mathcal{G} \coloneqq \trace(\nabla_{\cdot} \gradient_\mathcal{G}).
\end{equation}
\end{definition}
These Laplace operators clearly generalise the Laplace-Beltrami operator, and could be interesting on manifolds which have natural connections that are \emph{not} the Levi-Civita connection. In particular, on Lie groups there exists a family of canonical connections called the \emph{Lie-Cartan connections}.
\begin{definition}[Lie-Cartan Connection]\label{def:lie_cartan_connection}
Let $G$ be a Lie group and let $\nu \in \R$. Then the $\nu$ \emph{Lie-Cartan connection} $\nabla^{[\nu]}$ is the affine connection such that
\begin{equation}\label{eq:lie_cartan_connection}
\nabla^{[\nu]}_X Y = \nu [X, Y]
\end{equation}
for any left-invariant $X, Y \in \mathfrak{X}(G)$.
\end{definition}
Lie-Cartan connections have nice properties, such as the fact that their geodesics are exactly the exponential curves of $G$ \cite[Def.~2]{duits2021cartanconnection}. Of particular interest is $\nabla^{[0]}$, since this is the only one that is metric compatible with every left-invariant metric tensor field on $G$ \cite[Cor.~2]{duits2021cartanconnection}. We will now express the Laplace operators induced by Lie-Cartan connections and left-invariant metric tensor fields -- which we call \emph{Lie-Cartan Laplacians} -- and compare these with the Laplace-Beltrami operator. For readability, we write $\laplace_{\mathcal{G}, \nu} \coloneqq \laplace_{\mathcal{G}, \nabla^{[\nu]}}$. We also use the Einstein summation convention (e.g. \cite[App.~A]{lee2018riemannian}) in this section for concise  expressions. 
\begin{theorem}[Lie-Cartan Laplacians]\label{thm:left_invariant_lie_cartan_laplacian}
Let $G$ be a connected Lie group, let $\mathcal{G}$ be a left-invariant metric tensor field thereon, and let $\nu \in \R$. With respect to a left-invariant frame $\{\A_i\}_i$, the Lie-Cartan Laplacian is given by
\begin{equation}\label{eq:lie_cartan_laplacian_formula}
\laplace_{\mathcal{G}, \nu} = g^{ij} \A_i \A_j + \nu c_{ki}^k g^{ij} \A_j,
\end{equation}
while the Laplace-Beltrami operator is given by
\begin{equation}\label{eq:laplace_beltrami_formula}
\laplace_\mathcal{G} = g^{ij} \A_i \A_j + \Gamma_{ki}^k g^{ij} \A_j = g^{ij} \A_i \A_j + c_{ki}^k g^{ij} \A_j,
\end{equation}
with $c_{ij}^k$ the structure constants defined by $c_{ij}^k \A_k = [\A_i, \A_j]$. The difference is:
\begin{equation}\label{eq:difference_lie_cartan_laplacian_laplace_beltrami}
\laplace_\mathcal{G} - \laplace_{\mathcal{G}, \nu} = (1 - \nu) c_{ki}^k g^{ij} \A_j.
\end{equation}
\end{theorem}
\begin{proofsketch}
We first find for the Riemannian and Lie-Cartan divergences: 
\begin{equation*}
\begin{split}
\trace(\nabla_\cdot^\mathrm{LC} X) & = (\A_i + \Gamma_{ji}^j) X^i = (\A_i + c_{ji}^j) X^i, \textrm{ and} \\
\trace(\nabla_\cdot^{[\nu]} X) & = \A_i X^i + \nu X^j c_{ij}^j = (\A_i + \nu c_{ji}^j) X^i,
\end{split}
\end{equation*}
where $\Gamma_{ij}^k$ are the Christoffel symbols with respect to the left-invariant frame. Plugging in the gradient given by $\gradient_\mathcal{G} U = g^{ij} \A_i U \A_j$ gives the desired result.
\end{proofsketch}
Equation~\eqref{eq:difference_lie_cartan_laplacian_laplace_beltrami} tells us that the Laplace-Beltrami operator and the $1$ Lie-Cartan Laplacian coincide. The operators coincide for all $\nu \in \R$ if the underlying Lie group is \emph{unimodular}: for connected Lie groups this implies the trace of the structure constants $c_{ki}^k$ vanishes, allowing us to simplify Equations~\eqref{eq:lie_cartan_laplacian_formula} and \eqref{eq:laplace_beltrami_formula}: 
\begin{corollary}\label{cor:left_invariant_lie_cartan_laplacian}
Let $G$ be a connected unimodular Lie group and $\mathcal{G}$ a left-invariant metric tensor field thereon, and let $\nu \in \R$. Then the Lie-Cartan Laplacian and Laplace-Beltrami operator coincide:
\begin{equation}\label{eq:unimodular_difference_lie_cartan_laplacian_laplace_beltrami}
\laplace_{\mathcal{G}, \nu} = g^{ij} \A_i \A_j = \laplace_\mathcal{G}.
\end{equation}
\end{corollary}
Many Lie groups -- including $\SE(2)$ and $\Rtwo$ -- are unimodular.\footnote{On $\Rtwo$, all structure constants vanish; on $\SE(2)$, the nonzero structure constants are $c_{13}^2 = -c_{31}^2 = - c_{23}^1 = c_{32}^1$, none of which contribute to the trace of the adjoint.} As such, Cor.~\ref{cor:left_invariant_lie_cartan_laplacian} carries over onto their principal homogeneous spaces $\M$ and $\Rtwo$. 

The theory of Lie-Cartan connections can be generalised for gauge frames:
\begin{definition}[Gauge Frame Lie-Cartan Connection]\label{def:data_driven_lie_cartan_connection}
Let $G$ be a Lie group, with gauge frame $\{\A_i^U\}_i$, and let $\nu \in \R$. Then the $\nu$ \emph{gauge frame Lie-Cartan connection} $\nabla^{[\nu], U}$ is the affine connection such that
\begin{equation}\label{eq:data_driven_lie_cartan_connection}
\nabla^{[\nu], U}_X Y = \nu [X, Y]
\end{equation}
for any $X, Y \in \sections(T G)$ that are constant with respect to $\{\A_i^U\}_i$.
\end{definition}
We call the corresponding Laplacians \emph{data-driven Lie-Cartan Laplacians}, and write $\laplace_{\mathcal{G}, \nu}^U \coloneqq \laplace_{\mathcal{G}, \nabla^{[\nu], U}}$.
\begin{theorem}[Data-Driven Lie-Cartan Laplacians]\label{thm:data_driven_lie_cartan_laplacian}
Let $G$ be a Lie group, with gauge frame $\{\A_i^U\}_i$, let $\mathcal{G}$ be a  metric tensor field thereon that is constant with respect to the gauge frame, and let $\nu \in \R$. With respect to the gauge frame, the data-driven Lie-Cartan Laplacian is given by
\begin{equation}\label{eq:data_driven_lie_cartan_laplacian_formula}
\laplace_{\mathcal{G}, \nu}^U = g^{ij} \A_i^U \A_j^U + \nu d_{ki}^k g^{ij} \A_j^U,
\end{equation}
while the Laplace-Beltrami operator is given by
\begin{equation}\label{eq:data_driven_laplace_beltrami_formula}
\laplace_\mathcal{G} = g^{ij} \A_i^U \A_j^U + \Gamma_{ki}^k g^{ij} \A_j^U = g^{ij} \A_i^U \A_j^U + d_{ki}^k g^{ij} \A_j^U,
\end{equation}
with $d_{ij}^k $ the structure functions defined by $d_{ij}^k \A_k^U = [\A_i^U, \A_j^U]$. The difference is: 
\begin{equation}\label{eq:difference_data_driven_lie_cartan_laplacian_laplace_beltrami}
\laplace_\mathcal{G} - \laplace_{\mathcal{G}, \nu}^U = (1 - \nu) d_{ki}^k g^{ij} \A_j^U,
\end{equation}
so the two coincide if and only if $\nu = 1$.\footnote{In the gauge frame case, structure \emph{constants} $c_{ij}^k$ are replaced by structure \emph{functions} $d_{ij}^k$ which do not have a closed form and have non-vanishing trace in general.}
\end{theorem}
The proof is completely analogous to that of Theorem~\ref{thm:left_invariant_lie_cartan_laplacian}.

Hence, the diffusion we use for gauge frame DS filtering is not generated by the Laplace-Beltrami operator but by the $0$ data-driven Lie-Cartan Laplacian:
\begin{equation*}
\laplace \coloneqq g^{11} (\A_1^U)^2 + g^{22} (\A_2^U)^2 + g^{33} (\A_3^U)^2 \equiv \laplace_{\mathcal{G}, 0}^U.
\end{equation*}
One reason for this choice is that $\nabla^{[0], U}$ is the only gauge frame Lie-Cartan connection that is metric compatible with any metric tensor field that is constant with respect to the gauge frame. 
It moreover has a number of computational advantages compared to the Laplace-Beltrami operator: (1) it does not depend on the structure functions, and (2) all derivatives are of second order, whereas the Laplace-Beltrami operator has both second and first order derivatives. 
This data-driven Laplacian has been used in the past, e.g. \cite{franken2009cedos}. Similarly, the gauge frame TR-TV flow from \cite{smets2021tvflow} implicitly uses a divergence induced by $\nabla^{[0], U}$. We are interested in investigating the properties (e.g. analyticity) of the evolution generated by $\laplace_{\mathcal{G}, 0}^U$ in future work; in practice we see that the evolution is well-posed and smoothing as one would expect from a diffusion.
\section{Experimental Results}\label{sec:experiments}
We show a selection of experiments performed using RDS filtering on $\M$. Our Python implementations, which use Taichi \cite{taichi} for GPU acceleration, 
and experiments can be found at \url{https://github.com/finnsherry/M2RDSFiltering}.

\begin{figure}
\centering
\begin{tabular}{cc}
\includegraphics[width=0.45\linewidth]{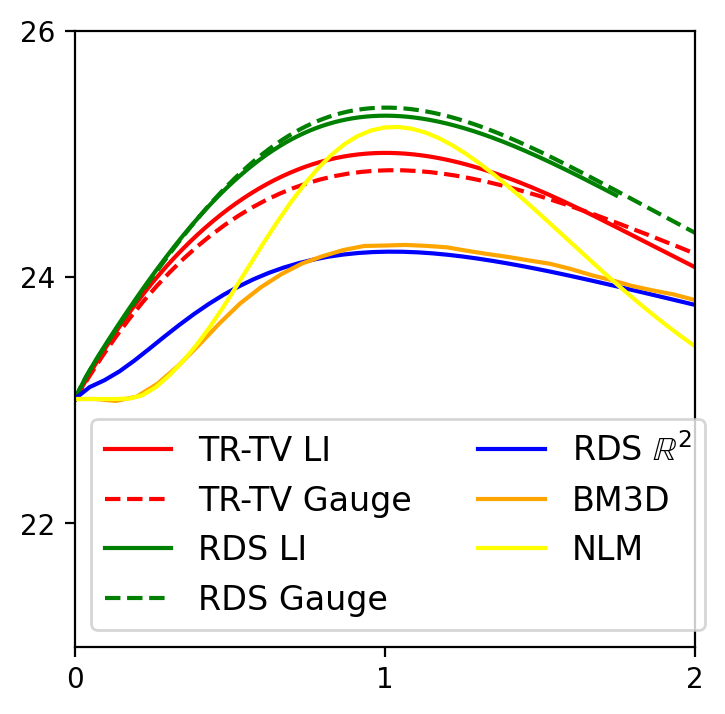}  &
\includegraphics[width=0.45\linewidth]{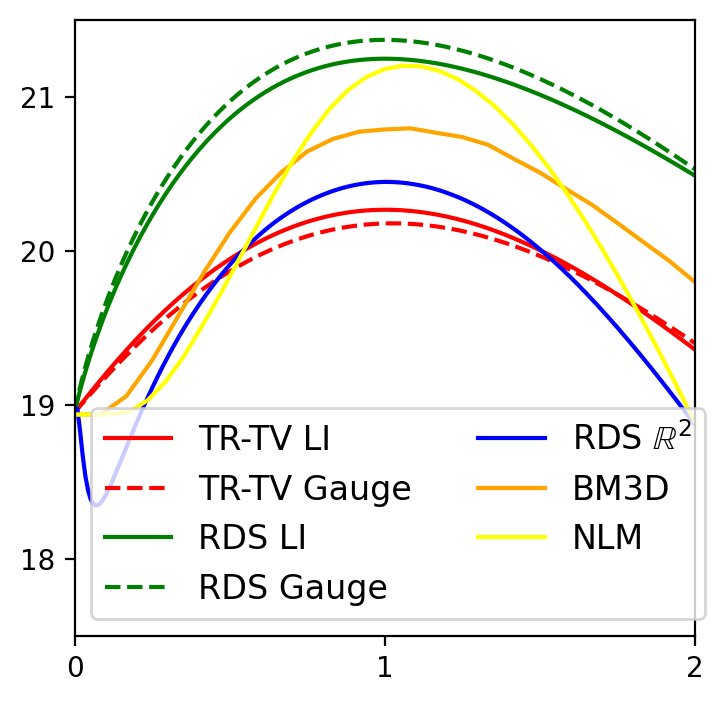} \\
(a) retina & (b) spiral
\end{tabular}
\caption{PSNRs of denoising methods over time relative to optimal stopping time.}
\label{fig:plots}
\end{figure}

\begin{figure}[tb]
\centering
\begin{tabular}{ccc}
    \small{(a) original} & (b) degraded image & (c) LI TR-TV on $\M$\\
     \includegraphics[width=0.3\linewidth]{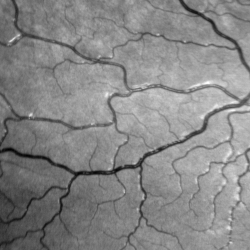} &
     \includegraphics[width=0.3\linewidth]{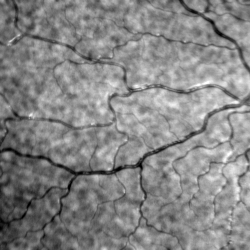}&
     \includegraphics[width=0.3\linewidth]{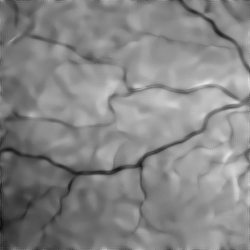}  \\
     \includegraphics[width=0.3\linewidth]{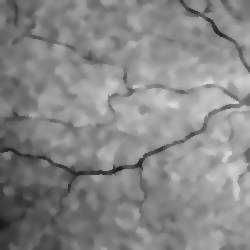}&
     \includegraphics[width=0.3\linewidth]{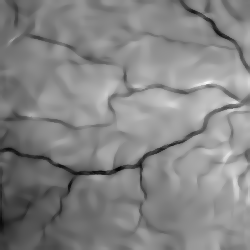}&
     \includegraphics[width=0.3\linewidth]{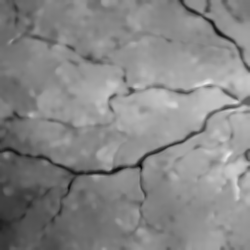}  \\
    (d) RDS on $\Rtwo$ & (e) Gauge RDS on $\M$& (f) NLM
\end{tabular}
\caption{Denoising of image \emph{retina} degraded with correlated noise.}
\label{fig:retina}
\end{figure}
\subsubsection{Image Enhancement.}
\begin{figure}[tb]
    \centering
    \begin{tabular}{ccc}
        \small{(a) original} & \small{(b) degraded image} & \small{(c) LI TR-TV on $\M$}\\
         \includegraphics[width=0.3\linewidth,frame]{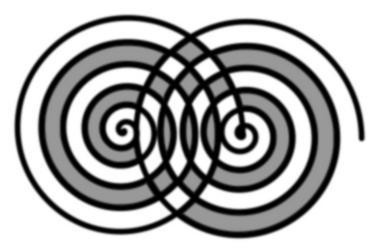} &
         \includegraphics[width=0.3\linewidth,frame]{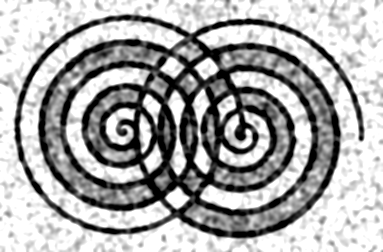}&
         \includegraphics[width=0.3\linewidth,frame]{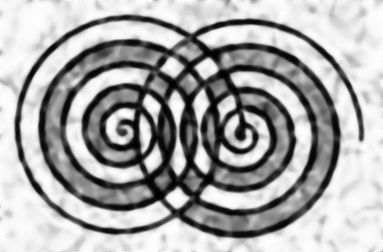}  \\
         \includegraphics[width=0.3\linewidth,frame]{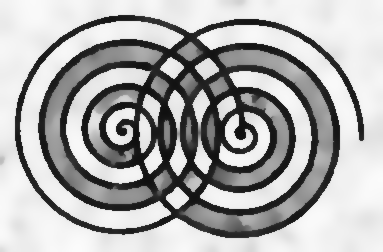}&
         \includegraphics[width=0.3\linewidth,frame]{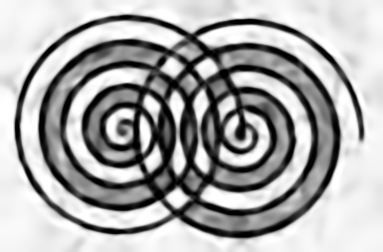}&
         \includegraphics[width=0.3\linewidth,frame]{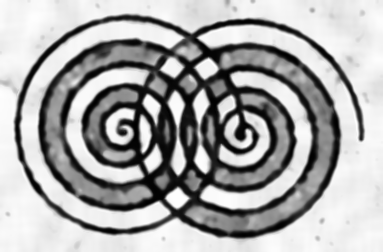}  \\
         \small{(d) RDS on $\Rtwo$} & \small{(c) Gauge RDS on $\M$}& \small{(d) NLM} 
    \end{tabular}
    \caption{Enhancement of image \emph{spiral} degraded by correlated noise.}
    \label{fig:spiral}
\end{figure}

We compare the enhancement of degraded images using gauge and left-invariant RDS filtering on $\M$ to our own implementations of RDS filtering on $\Rtwo$ \cite{schaefer2024regularisedds} and gauge and left-invariant TR-TV \cite{smets2021tvflow}; the \texttt{bm3d} Python implementation of BM3D \cite{bm3dold}; and the Python implementation of NLM \cite{nlm} in \texttt{scikit-image} \cite{scikit-image}. We consider a medical image \cite{Maastrichtstudy} -- a patch of the retina -- in Figure~\ref{fig:retina} and a cartoon-like image of overlapping spirals in Figure~\ref{fig:spiral}.  
In both cases the images have been degraded with additive, correlated noise $K_\rho * n_\sigma$, with $K_\rho$ a Gaussian of standard deviation $\rho$ and $n_\sigma$ white noise with intensity $\sigma$; in Figure~\ref{fig:retina} we have $(\sigma, \rho) = (127.5, 2)$ and in Figure~\ref{fig:spiral} we have $(\sigma, \rho) = (255, 2)$. 

Figure~\ref{fig:plots} shows the Peak Signal-to-Noise Ratio (PSNR) as a function of stopping time normalized by the optimal stopping time (for RDS, TR-TV) and as a function of noise power normalized to the optimal noise power (BM3D, NLM).
RDS filtering on $\M$ outperforms the other methods w.r.t. maximal PSNR, and is less sensitive to the stopping time than BM3D and NLM. 
In Figure~\ref{fig:retina} and \ref{fig:spiral} we qualitatively compare the results of various methods at their highest PSNRs.
We see in Figure~\ref{fig:retina} that the PDE-based methods preserve the small vessels better than NLM.
In Figure~\ref{fig:spiral}, we see that NLM has enhanced spurious blobs on the background which have been removed by $\M$ RDS.
Gauge frame RDS performs slightly better than left-invariant RDS.

\subsubsection{Inpainting.}
\begin{figure}[tb]
    \centering
    \begin{tabular}{ccc}
         \includegraphics[width=0.3\linewidth,frame]{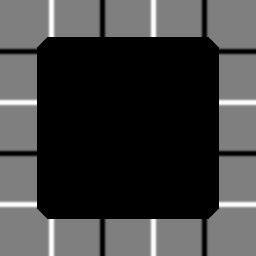} &
         \includegraphics[width=0.3\linewidth,frame]{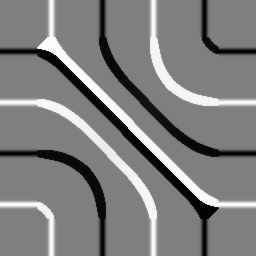}&
         \includegraphics[width=0.3\linewidth,frame]{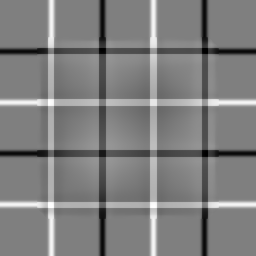}  \\
         (a) input & (b) $\Rtwo$ &  (c) $\M$\\
    \end{tabular}
    \caption{RDS inpainting of crossings on $\Rtwo$ and $\M$. The inpainting region is marked by a black square in (a).}
    \label{fig:inpainting}
\end{figure}

RDS filtering on $\Rtwo$ was originally developed with the goal of inpainting~\cite{BSCB00}, the task of filling in missing gaps in an image. 
While it can create good inpainting results, RDS inpainting on $\Rtwo$ -- like many other methods on $\Rtwo$ -- cannot reconstruct crossings. 
However, by lifting the image to $\M$ with the orientation score transform \eqref{eq:ost}, crossing structures are disentangled, allowing RDS inpainting on $\M$ to reconstruct crossings. 
Figure~\ref{fig:inpainting} demonstrates this: RDS inpainting Figure~\ref{fig:inpainting}c on $\M$ creates crossing lines, while on $\Rtwo$ it connects different lines without crossings Figure~\ref{fig:inpainting}b.
\subsubsection{Conclusion.}
We improved RDS filtering by extending it to $\M$, making it crossing-preserving. 
We created a left-invariant and gauge frame variant of RDS. We showed experimentally that both methods can outperform existing algorithms (TR-TV, BM3D, NLM, $\Rtwo$ RDS) in terms of maximal PSNR and sensitivity to stopping time, cf. Figure~\ref{fig:plots}. 
We showed that the Lie-Cartan and Levi-Civita connections induce the same Laplacian in the left-invariant situation (Theorem~\ref{thm:left_invariant_lie_cartan_laplacian}), but different ones in the gauge frame case (Theorem~\ref{thm:data_driven_lie_cartan_laplacian}).

\noindent
\textbf{Future work:} As crossing-preserving shock filtering is beneficial in geometric image processing on $\M$, we will integrate it in PDE-based deep learning, cf.~\cite{smets2024thesis}.

\subsubsection{Acknowledgements.}
We gratefully acknowledge the Dutch Foundation of Science NWO for its financial support by Talent Programme VICI 2020 Exact Sciences (Duits, Geometric learning for Image Analysis, VI.C 202-031).

\bibliographystyle{splncs04}
\bibliography{references}

\begin{thebibliography}{10}
\providecommand{\url}[1]{\texttt{#1}}
\providecommand{\urlprefix}{URL }
\providecommand{\doi}[1]{https://doi.org/#1}

\bibitem{AGLM93}
Alvarez, L., Guichard, F., Lions, P.L., Morel, J.M.: Axioms and fundamental
  equations in image processing. Arch. Ration. Mech. Anal.  \textbf{123},
  199--257 (1993). \doi{10.1007/BF00375127}

\bibitem{AVK93}
Arehart, A.B., Vincent, L., Kimia, B.B.: Mathematical morphology: The
  {H}amilton--{J}acobi connection. In: 4th Int. Conf. Comput. Vis. pp. 215--219
  (1993). \doi{10.1109/ICCV.1993.378217}

\bibitem{Prandi}
Baspinar, E., Calatroni, L., Franceschi, V., Prandi, D.: {A Cortical-Inspired
  Sub-Riemannian Model for Poggendorff-Type Visual Illusions}. Imaging
  \textbf{7} (2021). \doi{10.3390/jimaging7030041}

\bibitem{bekkers2015subriemanniangeodesics}
Bekkers, E.J., Duits, R., Mashtakov, A., Sanguinetti, G.R.: {A PDE Approach to
  Data-Driven Sub-Riemannian Geodesics in SE(2)}. SIIMS  \textbf{8},
  2740--2770 (2015). \doi{10.1137/15M1018460}

\bibitem{paperGijs}
Bellaard, G., Bon, D.L., Pai, G., Smets, B.M., Duits, R.: {Analysis of
  (sub-)Riemannian PDE-G-CNNs}. JMIV  \textbf{65},  819--843 (2023).
  \doi{10.1007/s10851-023-01147-w}

\bibitem{BSCB00}
Bertalm\'io, M., Sapiro, G., Caselles, V., Ballester, C.: Image inpainting. In:
  Proc. {SIGGRAPH} 2000. pp. 417--424. New Orleans (Jul 2000).
  \doi{10.1109/TIP.2003.815261}

\bibitem{BM92}
Brockett, R.W., Maragos, P.: Evolution equations for continuous-scale
  morphology. In: IEEE Int. Conf. Acoust. Speech Signal Process. vol.~3, pp.
  125--128 (1992). \doi{10.1109/78.340774}

\bibitem{nlm}
Buades, A., Coll, B., Morel, J.M.: {NLM Denoising}. IPOL pp. 208--212 (2011).
  \doi{10.5201/ipol.2011.bcm_nlm}

\bibitem{chambolle2019tv}
Chambolle, A., Pock, T.: {Total roto-translational variation}. Num. Math
  \textbf{142},  611--666 (2019). \doi{10.1007/s00211-019-01026-w}

\bibitem{charbonnier1997switch}
Charbonnier, P., Blanc-Féraud, L., Aubert, G., Barlaud, M.: {Deterministic
  Edge-Preserving Regularization in Computed Imaging}. TIP  \textbf{6},
  298--311 (1997). \doi{10.1109/83.551699}

\bibitem{cittisartiJMIV}
Citti, G., Franceschiello, B., Sanguinetti, G., Sarti, A.: {A Cortical Based
  Model of Perceptual Completion in the Roto-Translation Space}. JMIV
  \textbf{24},  307--326 (2006). \doi{10.1007/s10851-005-3630-2}

\bibitem{Citti2016mcf}
Citti, G., Franceschiello, B., Sanguinetti, G., Sarti, A.: {Sub-Riemannian Mean
  Curvature Flow for Image Processing}. SIIMS  \textbf{9},  212--237 (2016).
  \doi{10.1137/15M1013572}

\bibitem{bm3dold}
Dabov, K., Foi, A., Katkovnik, V., Egiazarian, K.: {mage Denoising by Sparse
  3-D Transform-Domain Collaborative Filtering}. TIP  \textbf{16},  2080--2095
  (2007). \doi{10.1109/TIP.2007.901238}

\bibitem{DuitsJMIV2013}
Duits, R., Dela~Haije, T., Creusen, E., Ghosh, A.: {Morphological and Linear
  Scale Spaces for Fiber Enhancement in DW-MRI}. JMIV  \textbf{46},  326--368
  (2013). \doi{10.1007/s10851-012-0387-2}

\bibitem{duits2021cartanconnection}
Duits, R., Smets, B.M., Wemmenhove, A., Portegies, J.W., Bekkers, E.J.: {Recent
  Geometric Flows in Multi-orientation Image Processing via a Cartan
  Connection}. Handb. of Math. Mod. Alg. Comp. Vis. Imaging pp. 1--60 (2021).
  \doi{10.1007/978-3-030-03009-4_101-1}

\bibitem{FG87}
F\"orstner, W., G\"ulch, E.: A fast operator for detection and precise location
  of distinct points, corners and centres of circular features. In: ISPRS
  Intercommiss. Conf. Fast Process. Photogrammetric Data. pp. 281--305 (1987)

\bibitem{franken2009cedos}
Franken, E.M., Duits, R.: {Crossing-Preserving Coherence-Enhancing Diffusion on
  Invertible Orientation Scores}. IJCV  \textbf{75},  253--278 (2009).
  \doi{10.1007/s11263-009-0213-5}

\bibitem{taichi}
Hu, Y., Li, T.M., Anderson, L., Ragan-Kelley, J., Durand, F.: Taichi: a
  language for high-performance computation on spatially sparse data
  structures. TOG  (2019)

\bibitem{Ii62}
Iijima, T.: Basic theory on normalization of pattern (in case of typical
  one-dimensional pattern). Bull. Electrotech. Lab.  \textbf{26},  368--388
  (1962)

\bibitem{lee2018riemannian}
Lee, J.M.: {Introduction to Riemannian Manifolds}. Springer (2018).
  \doi{10.1007/978-3-319-91755-9}

\bibitem{Liu2023ElasticaRegularization}
Liu, H., Tai, X.C., Kimmel, R., Glowinski, R.: {Elastica Models for Color Image
  Regularization}. SIIMS  \textbf{16}(1) (2023). \doi{10.1137/22M147935X}

\bibitem{petitotbook}
Petitot, J.: {Elements of Neurogeometry: Functional Architectures of Vision}.
  Springer (2017). \doi{10.1007/978-3-319-65591-8}

\bibitem{CalatroniSIAMreview2023}
Pragliola, M., Calatroni, L., Lanza, A., Sgallari, F.: {On and Beyond Total
  Variation Regularization in Imaging: The Role of Space Variance}. Review
  \textbf{65},  601--685 (2023). \doi{10.1137/21M1410683}

\bibitem{rouy1992viscosity}
Rouy, E., Tourin, A.: {A Viscosity Solution Approach to Shape-From-Shading}.
  SINUM  \textbf{29},  867--884 (1992). \doi{10.1137/0729053}

\bibitem{SW23}
Schaefer, K., Weickert, J.: Diffusion-shock inpainting. In: SSVM, Lect. Notes
  Comput. Sci., vol. 14009, pp. 588--600. Springer (2023).
  \doi{doi.org/10.1007/978-3-031-31975-4_45}

\bibitem{schaefer2024regularisedds}
Schaefer, K., Weickert, J.: {Regularised Diffusion-Shock Inpainting}. JMIV
  \textbf{66},  447--463 (2024). \doi{10.1007/s10851-024-01175-0}

\bibitem{Maastrichtstudy}
Schram, M., et~al.: {The Maastricht Study: an extensive phenotyping study on
  determinants of type 2 diabetes, its complications and its comorbidities}.
  European Journal of Epidimology  \textbf{29}(1),  439--451 (2014).
  \doi{10.1007/s10654-014-9889-0}

\bibitem{smets2024thesis}
Smets, B.M.: {Geometric Partial Differential Equations in Deep Learning and
  Image Processing}. Ph.D. thesis, TU/e (2024)

\bibitem{smets2021tvflow}
Smets, B.M., Portegies, J.W., St-Onge, E., Duits, R.: {Total Variation and Mean
  Curvature PDEs on $\mathbb{M}_d$}. JMIV  \textbf{63},  237--262 (2021).
  \doi{10.1007/s10851-020-00991-4}

\bibitem{So04}
Soille, P.: Morphological Image Analysis. Springer, 2 edn. (2004).
  \doi{10.1007/978-3-662-05088-0}

\bibitem{stonge2019mcf}
St-Onge, E., Meesters, S., Bekkers, E.J., Descoteaux, M., Duits, R.: {HARDI
  denoising with mean-curvature enhancement PDE on SE(3)}. In: 27th ISMRM
  (2019)

\bibitem{scikit-image}
van~der Walt, S., {S}ch\"onberger, J.L., {Nunez-Iglesias}, J., {B}oulogne, F.,
  {W}arner, J.D., {Y}ager, N., {G}ouillart, E., {Y}u, T., contributors:
  scikit-image. PeerJ  \textbf{2}, ~e453 (2014). \doi{10.7717/peerj.453}

\bibitem{We03}
Weickert, J.: Coherence-enhancing shock filters. In: Pattern Recognit., Lect.
  Notes Comput. Sci., vol.~2781, pp.~1--8. Springer (2003).
  \doi{10.1007/978-3-540-45243-0_1}

\end{thebibliography}
\end{document}